\documentclass[aop,seceqn,citesort,MSNbibl,dvips]{arximspdf}

\doi{10.1214/09-AOP505}
\volume{38}
\issue{4}
\pubyear{2010}
\firstpage{1492}
\lastpage{1506}

\makeatletter

\newtheorem{prop}{Proposition}[section]
\newtheorem{theorem}[prop]{Theorem}
\newtheorem{lem}[prop]{Lemma}
\newproclaim{rem}{Remark}
\def\bsuffix #1{#1}
\makeatother

\begin{document}
\begin{frontmatter}

\title{Local limit theorems in free probability theory}
\runtitle{Free local limit theorems}

\begin{aug}
\author[A]{\fnms{Jiun-Chau} \snm{Wang}\corref{}\thanksref{t1}\ead[label=e1]{jiuwang@mast.queensu.ca}}
\runauthor{J.-C. Wang}
\affiliation{Queen's University}
\address[A]{Department of Mathematics\\
\quad and Statistics\\
Queen's University\\
Kingston, Ontario K7L 3N6\\
Canada\\
\printead{e1}} 
\end{aug}

\thankstext{t1}{Supported in part by a Coleman Postdoctoral
Fellowship at Queen's University.}

\received{\smonth{5} \syear{2009}}
\revised{\smonth{8} \syear{2009}}

%
\begin{abstract}
In this paper, we study the superconvergence phenomenon in the free
central limit theorem for identically distributed, unbounded
summands. We prove not only the uniform convergence of the densities
to the semicircular density but also their $L^{p}$-convergence to
the same limit for $p>1/2$. Moreover, an entropic central limit
theorem is obtained as a consequence of the above results.
\end{abstract}

%

\end{frontmatter}

\section{Introduction}

Let $\{ X_{i}\}_{i=1}^{\infty}$ be a sequence of freely independent
copies of a selfadjoint random variable $X$. For notational
convenience, we assume that the random variable $X$ has mean zero
and variance one. Let $\mu_{n}$ be the distribution of the
normalized sum $S_{n}=(X_{1}+X_{2}+\cdots+X_{n})/\sqrt{n}$. The aim
of this article is to study the convergence properties of the
probability measures $\mu_{n}$ in terms of their densities
$d\mu_{n}/dx$.

Convergence of the densities of a sequence of distributions to the
density of the limit distribution (assuming the given densities
exist) is the object of classical local limit theorems \cite{GK}.
One such limit theorem concerning the sequence
$\{\mu_{n}\}_{n=1}^{\infty}$ was discovered by Bercovici and
Voiculescu in \cite{Superconvergence}, where the random variable $X$
is assumed to be \textit{bounded} (i.e., the law of $X$ has a bounded
support). The limit law in this case is the standard semicircular
distribution $\gamma$ whose density is $\sqrt{4-x^{2}}/2\pi$ on the
interval $[-2,2]$. This convergence to the semicircle law, referred
as the \textit{superconvergence} in \cite{Superconvergence}, is very
strong in the sense that the measures $\mu_{n}$ become Lebesgue
absolutely continuous after finitely many $n$'s, the density of
$\mu_{n}$ is supported on an interval $[a_{n},b_{n}]$ and analytic
on $(a_{n},b_{n})$, $\lim_{n\rightarrow\infty}a_{n}=-2$,
$\lim_{n\rightarrow\infty}b_{n}=2$, and the sequence $d\mu_{n}/dx$
converges uniformly to $d\gamma/dx$ as $n\rightarrow\infty$.

On the other hand, the weak convergence of the sequence
$\{\mu_{n}\}_{n=1}^{\infty}$ to the measure $\gamma$, also known as
the \textit{free central limit theorem}, has been the subject of
several investigations. This limit theorem was initially proved by
Voiculescu \cite{VoiSym} for bounded random variables. Later this
result was extended to unbounded variables with finite variance by
Maassen \cite{Maassen} (see also \cite{PataCLT} for a further
development). An explanation of the appearance of the semicircle
law, both in the free central limit theorem and in the asymptotics
of large random matrices, was found by Voiculescu in \cite{VoiRM}.

Knowing the above results, it is therefore natural to ask whether
the superconvergence in the free central limit theorem holds for
unbounded random variables. In this paper, we show that this
is indeed the case and even more, in addition to the uniform
convergence and analyticity properties of the densities, we also
prove a global central limit theorem in the sense of the
$L^{p}$-convergence. This further implies an entropic central limit
theorem, which is a free analogue of Barron's classical theorem on
the convergence to the normal entropy \cite{Barron}. Note that our
results are in sharp contrast with the classical central limit
theorem, where the distributions in the central limit process could
be all atomic.

This paper is written in terms of free convolution. Let
$\mathcal{M}$ be the set of all Borel probability measures on the
real line $\mathbb{R}$. The free convolution of two measures
$\mu,\nu\in\mathcal{M}$ is denoted by $\mu\boxplus\nu$. Thus,
$\mu\boxplus\nu$ is the probability distribution of $Y+Z$, where $Y$
and $Z$ are free random variables with distributions $\mu$ and
$\nu$, respectively (see \cite{BVUnbdd} for the details of this
construction and the book \cite{DVN} for a comprehensive
introduction to free probability theory).

Aside from being an interesting subject in itself, the theory of
free convolution also plays a significant role in the study of the
limit laws for large random matrices. For example, it was shown in
\cite{PasturandVasilchuk} that if for each $N\geq1$, $A_{N}$ and
$B_{N}$ are independent unitarily invariant $N\times N$ hermitian
random matrices such that their spectral distributions converge
weakly in probability to $\mu$ and $\nu$ as $N\rightarrow\infty$ and
the sequence $\mathrm{E}(\operatorname{tr} |A_{N} |)$ is bounded
(where $\mathrm{E}$ denotes the expectation and $\operatorname{tr}$ the
normalized trace), then the spectral distribution of $A_{N}+B_{N}$
converges weakly in probability to the free convolution
$\mu\boxplus\nu$ as $N\rightarrow\infty$. Hence, one can use the
free probability tools such as the $R$-transform (see Section \ref{sec2}) to
further analyze the limit law of $A_{N}+B_{N}$. We refer to
\cite{BianeSurvey} for an introduction to some of the most
probabilistic aspects of free probability and its connections with
random matrices. We now begin by reviewing the analytical machinery
needed for the calculation of free convolution.

\section{Preliminaries}\label{sec2}

\subsection{Cauchy transforms and free convolution}

Denote by $\mathbb{C}^{+}=\{ z\in\mathbb{C}\dvtx\Im z>0\}$ the complex
upper half plane. The \textit{Cauchy transform} of a measure
$\mu\in\mathcal{M}$ is
\[
G_{\mu}(z)=\int_{-\infty}^{\infty}\frac{1}{z-x} \, d\mu(x),\qquad
z\in\mathbb{C}^{+},
\]
and its reciprocal $F_{\mu}=1/G_{\mu}$ is an
analytic self-mapping of $\mathbb{C}^{+}$. The measure $\mu$ can be
recovered from $G_{\mu}$ as the $\mathrm{weak}^{*}$-limit of the
measures
\[
d\mu_{y}(x)=-\frac{1}{\pi}\Im G_{\mu}(x+iy) \, dx
\]
as
$y\rightarrow0^{+}$ (see \cite{Donoghue}). If the function $\Im
G_{\mu}(z)$ is continuous at $x\in\mathbb{R}$, then the probability
distribution function $\mathcal{F}_{\mu}(t)=\mu((-\infty,t])$ is
differentiable at $x$ and its derivative is given by
$\mathcal{F}_{\mu}^{\prime}(x)=-\Im G_{\mu}(x)/\pi$. This inversion
formula gives a way to extract the density function of a measure
from its Cauchy transform.

The Cauchy transform $G_{\mu}$ carries certain invertibility
properties near the point of infinity, which can be used to
calculate the free convolution of measures. To be precise, for each
$\theta\in(0,\pi/2)$ and $\beta>0$, we introduce the cone
$\Gamma_{\theta}=\{ z\in\mathbb{C}^{+}\dvtx\arg
z\in(\theta,\pi-\theta)\}$ and the domain
\[
D_{\theta,\beta}=\{
z\in\mathbb{C}^{-}\dvtx\arg
z\in(-\pi+\theta,-\theta);  |z |<\beta\},
\]
where $\arg z$
is the principal argument of the complex number $z$ and
$\mathbb{C}^{-}$ denotes the complex lower half plane. As shown in
\cite{BVUnbdd}, we have
\[
G_{\mu}(z)=\frac{1}{z}\bigl(1+o(1)\bigr),\qquad
z\in\mathbb{C}^{+},
\]
as $z\rightarrow\infty$ \textit{nontangentially}
(i.e., $ |z |\rightarrow\infty$ but $z\in\Gamma_{\theta}$
for some angle $\theta>0$). This implies that, for every $\theta>0$,
there exists $\beta=\beta(\mu,\theta)>0$ such that the function
$G_{\mu}$ has an analytic inverse $G_{\mu}^{-1}$ (relative to the
composition) defined in the set~$D_{\theta,\beta}$. Then the
\textit{$R$-transform} of the measure $\mu$ is defined as
\[
R_{\mu}(z)=G_{\mu}^{-1}(z)-\frac{1}{z}
\]
and we have
\[
R_{\mu_{1}\boxplus\mu_{2}}=R_{\mu_{1}}+R_{\mu_{2}}
\]
on a domain
$D_{\theta,\beta}$ where these three functions are defined. This
remarkable property of $R$-transform was first proved by Voiculescu
in the case of compactly supported measures \cite{Voiadd}, then it
was extended to the measures with finite variance in~\cite{Maassen},
and finally to the whole class $\mathcal{M}$ \cite{BVUnbdd}.

Note that the $R$-transform $R_{\mu}$ is related with the function
$\varphi_{\mu}$ in \cite{BVUnbdd} via the formula
$\varphi_{\mu}(z)=R_{\mu} (1/z )$. Thus, the properties of
the function $\varphi_{\mu}$ can be translated in terms of
$R_{\mu}$. Among all, we mention the following continuity property
with respect to the weak convergence of probability measures (see
\cite{BVUnbdd}). The notation $\delta_{x}$ here means the
probability measure concentrated at $x\in\mathbb{R}$.
\begin{prop}\label{prop21}
Let $\{\nu_{n}\}_{n=1}^{\infty}\subset\mathcal{M}$ be a sequence. If
$\nu_{n}$ converges weakly to $\delta_{0}$, then for every
$\theta\in(0,\pi/2)$ and $\beta>0$ there exists
$N=N(\theta,\beta)>0$ such that the function $R_{\nu_{n}}$ is
defined in the domain $D_{\theta,\beta}$ for $n\geq N$ and
$R_{\nu_{n}}\rightarrow0$ uniformly on the compact subsets of
$D_{\theta,\beta}$ as $n\rightarrow\infty$.
\end{prop}

\subsection{Free convolution with a semicircle law}\label{sec22}

For $t>0$, the centered semicircle distribution of variance $t$ is
the probability measure with density
\[
d\gamma_{t}(x)=\frac{1}{2\pi
t}\sqrt{4t-x^{2}} \, dx
\]
on the interval $[-2\sqrt{t},2\sqrt{t}]$.
Its Cauchy transform is given by
\[
G_{\gamma_{t}}(z)=\frac{z-\sqrt{z^{2}-4t}}{2t},\qquad
z\in\mathbb{C}^{+},
\]
where the branch of the square root on
$\mathbb{C}^{+}\setminus[0,+\infty)$ is chosen such that
$\sqrt{-1}=i$. We use $\gamma$ to denote the standard semicircle law
that has mean zero and variance one. The function $G_{\gamma}$ has a
continuous extension to $\mathbb{C}^{+}\cup\mathbb{R}$, where the
extension acts on $\mathbb{R}$ by
\[
\cases{
 \bigl(x-i\sqrt{4-x^{2}} \bigr)/2, &\quad if $ |x |\leq2$;\cr
 \bigl(x-\sqrt{x^{2}-4} \bigr)/2, &\quad
if $ |x |>2$.}
\]
In particular, we see
that, for each $\delta>0$, the function $G_{\gamma}$ can be
continued analytically to $K=\{
x+iy\dvtx x\in(-2,2), |y |<\delta\}$ (and beyond, to a whole
Riemann surface). This analytic continuation (which we still denote
by $G_{\gamma}$) has an explicit formula given by
$G_{\gamma}(z)= (z-i(4-z^{2})^{1/2} )/2$, where
$(\cdot)^{1/2}$ is the principal branch of the square root function
on $\mathbb{C}^{+}\setminus(-\infty,0]$. Note that the function
$G_{\gamma_{t}}$ also has a similar boundary behavior since
$G_{\gamma_{t}}(z)=t^{-1/2}G_{\gamma} (t^{-1/2}z )$. We
also mention that the $R$-transform $R_{\gamma_{t}}(z)=tz$ and the
function $G_{\gamma}$ satisfies the functional
equation
%
%
\begin{equation}\label{eq:2.1}
G_{\gamma}(z)+F_{\gamma}(z)=z,\qquad
z\in\mathbb{C}^{+}\cup K.
\end{equation}

Proprieties of free convolution by semicircular distributions have
been studied thoroughly in \cite{VoiI,Biane}. We now review some of
these results that are relevant to our approach to the central limit
theorem. Fix $t>0$ and a measure $\nu\in\mathcal{M}$. As shown by
Biane \cite{Biane}, the Cauchy transform $G_{\nu\boxplus\gamma_{t}}$
has a continuous extension to $\mathbb{C}^{+}\cup\mathbb{R}$ and an
explicit formula for the density of the measure
$\nu\boxplus\gamma_{t}$ can be described as follows.

Define the function $v_{t}\dvtx\mathbb{R}\rightarrow[0,+\infty)$ by
\[
v_{t}(u)=\inf \biggl\{
v\geq0\dvtx \int_{-\infty}^{\infty}\frac{1}{(u-x)^{2}+v^{2}}
\,d\nu(x)\leq\frac{1}{t} \biggr\}.
\]
It was proved in \cite{Biane}, Lemma 2, that the function $v_{t}$ is
continuous on $\mathbb{R}$,
analytic on the open set $\{ u\in\mathbb{R}\dvtx v_{t}(u)>0\}$ and
for all $u\in\mathbb{R}$
%
%
\begin{equation}\label{eq:2.2}
\int_{-\infty}^{\infty}\frac{1}{(u-x)^{2}+v_{t}(u)^{2}}\,
d\nu(x)\leq\frac{1}{t}
\end{equation}
with equality if
$v_{t}(u)>0$. Let
\[
\psi_{t}(u)=u+t\int_{-\infty}^{\infty}\frac
{(u-x)}{(u-x)^{2}+v_{t}(u)^{2}}\,
d\nu(x), \qquad u\in\mathbb{R}.
\]
The following formula is proved by
Biane \cite{Biane}.
\begin{theorem}\label{theo22}
The function $\psi_{t}\dvtx\mathbb{R}\rightarrow\mathbb{R}$ is a
homeomorphism and at the point $x=\psi_{t}(u)$ we have
\[
G_{\nu\boxplus\gamma_{t}}(x)=\frac{1}{t} [x-u-iv_{t}(u) ].
\]
\end{theorem}

Thus, the density of the measure $\nu\boxplus\gamma_{t}$ is simply
$v_{t}(u)/\pi t$ at the point $\psi_{t}(u)$. Biane also showed that
the function $G_{\nu\boxplus\gamma_{t}}(z)$ has an analytic
extension to wherever $v_{t}$ is positive (see Corollary 4 in
\cite{Biane}).

Finally, we note for a further reference the following estimate in
\cite{Biane}:
%
%
\begin{equation}\label{eq:2.3}
|G_{\nu\boxplus\gamma_{t}}(z) |\leq\frac{1}{\sqrt
{t}},\qquad
z\in\mathbb{C}^{+}\cup\mathbb{R}.
\end{equation}

\subsection{Analytic subordination for free convolution powers}

For a measure $\mu$ in $\mathcal{M}$ and a positive integer
$n\geq2$, the $n$th free convolution power of $\mu$ is defined as
\[
\mu^{\boxplus
n}=\underbrace{\mu\boxplus\mu\boxplus\cdots\boxplus\mu}_{n\
\mathrm{times}}.
\]
As shown in \cite{Semigroup}, the reciprocal Cauchy transform
$F_{\mu^{\boxplus n}}$ is subordinated to the function $F_{\mu}$,
that is,
%
%
\begin{equation}\label{eq:2.4}
F_{\mu^{\boxplus
n}}(z)=F_{\mu} (\omega_{n}(z) ),\qquad
z\in\mathbb{C}^{+},
\end{equation}
where the
subordination function $\omega_{n}$ is given by the
formula
%
%
\begin{equation}\label{eq:2.5}
\omega_{n}(z)=\frac{z}{n}+\frac{n-1}{n}F_{\mu^{\boxplus
n}}(z).
\end{equation}

It was also proved in \cite{Semigroup} that the function
$F_{\mu^{\boxplus n}}$ has a continuous extension to
$\mathbb{C}^{+}\cup\mathbb{R}$ and the singular part in the Lebesgue
decomposition of the measure $\mu^{\boxplus n}$ is purely atomic.
Moreover, the measure $\mu^{\boxplus n}$ has at most one atom and a
number $\alpha_{n}$ is an atom of $\mu^{\boxplus n}$ if and only if
$x=\alpha_{n}/n$ is an atom of $\mu$ such that $\mu^{\boxplus
n}(\{\alpha_{n}\})=n\mu(\{ x\})-n+1$. Thus, if the measure $\mu$ is
not a point mass, then the atom of $\mu^{\boxplus n}$ disappears
when $n$ is sufficiently large and hence $\mu^{\boxplus n}$ is
Lebesgue absolutely continuous.

\section{Convergence of densities}\label{sec3}

Let $\mu$ be any probability measure on $\mathbb{R}$ with finite
nonzero variance. For simplicity's sake, we confine attention to the
case of zero mean and unit variance, that is,
\[
\int_{-\infty}^{\infty}x\,
d\mu(x)=0 \quad\mbox{and}\quad \int_{-\infty}^{\infty}x^{2}\,
d\mu(x)=1.
\]
Indeed, the general case can be reduced to this special
case by a simple translating and scaling argument. For such a
measure $\mu$ and a positive integer $n\geq2$, we introduce
\[
\mu_{n}=\underbrace{D_{{1}/{\sqrt{n}}}\mu\boxplus
D_{{1}/{\sqrt{n}}}\mu\boxplus\cdots\boxplus
D_{{1}/{\sqrt{n}}}\mu}_{n\ \mathrm{times}},
\]
where the measure
$D_{1/\sqrt{n}}\mu$ is the \textit{dilation} of $\mu$ by a factor of
$1/\sqrt{n}$, that is,
\[
dD_{{1}/{\sqrt{n}}}\mu(x)=d\mu \bigl(\sqrt{n}x \bigr).
\]
Note
that the measure $\mu_{n}$ can also be realized as the distribution
of
\[
\frac{1}{\sqrt{n}} (X_{1}+X_{2}+\cdots+X_{n} ),
\]
where $X_{1},X_{2},\ldots,X_{n}$ are freely independent random
variables affiliated to a $W^{*}$-probability space, all having
distribution $\mu$ (see \cite{BVUnbdd}). In other words, we have
$\mu_{n}=D_{1/\sqrt{n}}\mu^{\boxplus n}$. At the level of Cauchy
transforms, this means that $G_{\mu_{n}}(z)=\sqrt{n}G_{\mu^{\boxplus
n}} (\sqrt{n}z )$. Meanwhile, by Segal's noncommutative
integration theory \cite{Segal} and the work of Maassen
\cite{Maassen}, each measure $\mu_{n}$ also has mean zero and
variance one. The free central limit theorem \cite{Maassen} states
that the sequence $\mu_{n}$ converges weakly to $\gamma$ as
$n\rightarrow\infty$. In this section, we will show that the density
functions $d\mu_{n}/dx$ actually converge to $d\gamma/dx$ uniformly
on $\mathbb{R}$ as well as in $L^{p}$-norms for $p>1/2$.

To our purposes, we first reinterpret the central limit theorem in
terms of free Brownian motion. The following result is essentially
the same as Theorem 1.6 in~\cite{SerNica}. We provide its proof here
for the sake of completeness. Throughout this paper, we will use the
following notation:
\[
t=t(n)=\frac{n}{n+1} \quad\mbox{and}\quad
F_{n}(z)=F_{\mu_{n}}(z).
\]
\begin{lem}\label{lem31}
There exists a unique probability measure $\nu\in\mathcal{M}$ such
that
%
%
\begin{equation}\label{eq:3.1} F_{\mu}(z)=z-G_{\nu}(z),\qquad
z\in\mathbb{C}^{+},
\end{equation}
and, for every
$n\geq1$,
%
%
\begin{equation}\label{eq:3.2}
F_{n+1}(z)=z-G_{\nu_{n}\boxplus\gamma_{t}}(z),\qquad
z\in\mathbb{C}^{+}\cup\mathbb{R},
\end{equation}
where
the measure $\nu_{n}$ is given by
$d\nu_{n}(x)=d\nu (\sqrt{n+1}x )$.
\end{lem}
\begin{pf}
The proof of the representation (\ref{eq:3.1}) can be found in
\cite{Maassen}, Proposition~2.2.

We will focus on the proof of (\ref{eq:3.2}). Notice that we only
need to show that (\ref{eq:3.2}) holds in an open subset of
$\mathbb{C}^{+}$ since both functions $F_{n+1}$ and
$G_{\nu_{n}\boxplus\gamma_{t}}$ extend continuously to
$\mathbb{C}^{+}\cup\mathbb{R}$, and they are analytic in
$\mathbb{C}^{+}$. Let $\omega_{n+1}(z)$ be the subordination
function of $F_{\mu^{\boxplus(n+1)}}(z)$ as in (\ref{eq:2.4}). Then
(\ref{eq:2.5}) and (\ref{eq:3.1}) imply that
\begin{eqnarray*}
\omega_{n+1}(z) & = & \frac{1}{n+1}z+\frac{n}{n+1}F_{\mu}
(\omega
_{n+1}(z) )\\
& = & \frac{1}{n+1}z+\frac{n}{n+1}\omega_{n+1}(z)-\frac
{n}{n+1}G_{\nu
} (\omega_{n+1}(z) ).
\end{eqnarray*}
Consequently, we have
%
%
\begin{equation}\label{eq:3.3}
\omega_{n+1}(z)+nG_{\nu} (\omega_{n+1}(z) )=z,\qquad
z\in\mathbb{C}^{+}.
\end{equation}

On the other hand, the additivity of $R$-transform shows that there
exists a large $\beta>0$ such that the equation
\[
G_{\nu\boxplus\gamma_{n}} \bigl(w+nG_{\nu}(w) \bigr)=G_{\nu}(w)
\]
holds in a truncated cone $\Gamma=\{
x+iy\in\mathbb{C}^{+}\dvtx |x |<y;  y>\beta\}$. Combining this
with (\ref{eq:3.3}), we obtain the identity
%
%
\begin{equation}\label{eq:3.4}
G_{\nu\boxplus\gamma_{n}} (z )=G_{\nu} (\omega
_{n+1}(z) )
\end{equation}
in a neighborhood of $\infty$; notice that we have used implicitly a
consequence of~(\ref{eq:2.5}), namely, the function
$\omega_{n+1}(z)=z(1+o(1))$ as $z\rightarrow\infty$ nontangentially
so that $\omega_{n+1}(z)$ lies in $\Gamma$ as $z\rightarrow\infty$,
$z\in\Gamma$.

Meanwhile, (\ref{eq:2.4}) and (\ref{eq:3.1}) show that
\[
F_{\mu^{\boxplus(n+1)}}(z)=F_{\mu} (\omega_{n+1}(z)
)=\omega
_{n+1}(z)-G_{\nu} (\omega_{n+1}(z) ),\qquad
z\in\mathbb{C}^{+}.
\]
Together with (\ref{eq:3.3}) and
(\ref{eq:3.4}), we deduce that the identity
\[
F_{\mu^{\boxplus(n+1)}}(z)=z-(n+1)G_{\nu\boxplus\gamma_{n}}
(z )
\]
holds in a neighborhood of $\infty$. After dilating by a factor of
$1/\sqrt{n+1}$, we conclude that (\ref{eq:3.2}) holds in an open
subset of $\mathbb{C}^{+}$ and hence the proof is complete.
\end{pf}

Following Biane \cite{Biane}, we now parametrize the real line by
the homeomorphism $\psi_{t}$ as in Section \ref{sec22}. Setting
$x=\psi_{t}(u)$, we first derive an estimate for the range of~$x$.
\begin{lem}\label{lem32}
For every $\eta\in(0,1)$ and $n\geq1$, the set
$\psi_{t} ([-1+\eta,1-\eta] )$ is contained in the interval
$[-2+\eta,2-\eta]$.
\end{lem}
\begin{pf}
By Theorem \ref{theo22}, we have $\psi_{t}(u)=u+t\Re
G_{\nu_{n}\boxplus\gamma_{t}}(x)$. Moreover, (\ref{eq:2.3}) implies
that
\[
 |t\Re
G_{\nu_{n}\boxplus\gamma_{t}}(x) |\leq\sqrt{t}=\sqrt{\frac
{n}{n+1}}<1-\frac{1}{2(n+1)}.
\]
Hence, we obtain
\[
 |\psi_{t}(u) |<2-\eta-\frac{1}{2(n+1)},\qquad
u\in[-1+\eta,1-\eta],
\]
which implies the desired result.
\end{pf}

We need one more detail about the boundary behavior of the function
$G_{\nu_{n}\boxplus\gamma_{t}}(z)$ as $n\rightarrow\infty$. Observe
that the sequence $\nu_{n}$ converges weakly to $\delta_{0}$.
\begin{lem}\label{lem33}
For each $\eta\in(0,1)$ there exists $N=N(\eta)>0$ such that, for
all $n\geq N$, the function $G_{\nu_{n}\boxplus\gamma_{t}}(z)$ can
be continued analytically to a neighborhood of the interval
$[-2+\eta,2-\eta]$. Furthermore, this analytic continuation never
vanishes on $[-2+\eta,2-\eta]$.
\end{lem}
\begin{pf}
By virtue of Theorem \ref{theo22}, it suffices to show that $v_{t}(u)>0$ on
the interval $[-1+\eta,1-\eta]$ for sufficiently large $n$. Suppose
that, on the contrary, we can find a subsequence $\{
n_{k}\}_{k=1}^{\infty}\subset\mathbb{N}$ and a sequence $\{
u_{k}\}_{k=1}^{\infty}$ in $[-1+\eta,1-\eta]$ such that
$n_{1}<n_{2}<\cdots$ and $v_{t_{k}}(u_{k})=0$ for all $k$'s, where
$t_{k}$ denotes $n_{k}/(n_{k}+1)$. Then the weak convergence of
$\{\nu_{n}\}_{n=1}^{\infty}$ implies that
\[
\nu_{n_{k}} \biggl( \biggl(-\frac{\eta}{2},\frac{\eta}{2}
\biggr)
\biggr)>1-\frac{\eta}{2}
\]
for sufficiently large $k$. Therefore, we deduce, from
(\ref{eq:2.2}), that
\begin{eqnarray*}
\frac{1}{t_{k}} & > & \int_{-\infty}^{\infty}\frac
{1}{(u_{k}-x)^{2}}
\,d\nu_{n_{k}}(x)\\
& \geq& \int_{ |x |<{\eta}/{2}}\frac
{1}{(u_{k}-x)^{2}}
\,d\nu_{n_{k}}(x)>\frac{2}{2-\eta}
\end{eqnarray*}
for large $k$. This contradicts to the fact that
$\lim_{k\rightarrow\infty}t_{k}=1$. Hence, the function $v_{t}$ must
be positive on $[-1+\eta,1-\eta]$ when $n$ is large. 
\end{pf}

We now have all ingredients to proceed with the proof of our main
result.
\begin{theorem}\label{theo34}
Suppose that $\mu\in\mathcal{M}$ has mean zero and variance one.
Then:
\begin{enumerate}[(iii)]
\item[(i)] the measure $\mu_{n}$ is Lebesgue absolutely
continuous for sufficiently
large $n$;
\item[(ii)] for each small $\varepsilon>0$ there exist
$\delta>0$ and $N>0$
such that the function $G_{\mu_{n}}$ has an analytic continuation
$h_{n}$ to $K=\{
x+iy\dvtx x\in[-2+\varepsilon,2-\varepsilon], |y |<\delta\}$
whenever $n\geq N$. Moreover,
$h_{n}(z)\rightarrow(z-i(4-z^{2})^{1/2})/2$ uniformly on $K$ as
$n\rightarrow\infty$;
\item[(iii)] the density $d\mu_{n}/dx$ is continuous for
sufficiently large $n$
and $d\mu_{n}/dx\rightarrow d\gamma/dx$ uniformly on $\mathbb{R}$ as
$n\rightarrow\infty$.
\end{enumerate}
\end{theorem}
\begin{pf}
The statement (i) follows from the regularity results of the measure
$\mu^{\boxplus n}$ because $\mu_{n}$ is a dilation of $\mu^{\boxplus
n}$.

Let us prove (ii). In accordance with the established terminology,
we will work with measures indexed by $n+1$ instead of $n$. Let
$\varepsilon\in(0,1/10)$ be given. Define the sets
$K_{\varepsilon}=\{
x+iy\dvtx x\in[-2+\varepsilon,2-\varepsilon], |4y
|<\varepsilon\sqrt
{\varepsilon}\}$
and $D_{\theta,2}=\{ z\in\mathbb{C}^{-}\dvtx\arg
z\in(-\pi+\theta,-\theta);  |z |<2\}$, where the angle
$\theta=\theta(\varepsilon)$ is chosen so that
$4\sin\theta=\sqrt{\varepsilon(4-\varepsilon)}$.

We first show that the set $G_{\gamma} (K_{\varepsilon} )$
is contained in $D_{\theta,2}$, where
$G_{\gamma}(z)=(z-i(4-z^{2})^{1/2})/2$ is the analytic extension of
the Cauchy transform of $\gamma$ on $K_{\varepsilon}$. Let $z_{0}\in
K_{\varepsilon}$ be given. We write $G_{\gamma}(z_{0})=Re^{i\psi}$,
where $\psi=\arg G_{\gamma}(z_{0})$. To prove that
$G_{\gamma}(z_{0})$ is in $D_{\theta,2}$, we need to verify that
$ |\sin\psi |>\sin\theta$ and $R<2$. Indeed, from the
functional equation (\ref{eq:2.1}), we have
\[
\biggl(R+\frac{1}{R} \biggr)\cos\psi+i \biggl(R-\frac{1}{R}
\biggr)\sin\psi=z_{0}.
\]
Using the fact that $ |\Re z_{0} |\leq2-\varepsilon$, we
get
\[
2 |{\cos}\psi |\leq \biggl(R+\frac{1}{R} \biggr) |{\cos}
\psi
|\leq2-\varepsilon.
\]
This implies that
$2 |{\sin}\psi |\geq\sqrt{\varepsilon(4-\varepsilon
)}>2\sin\theta$,
as desired. On the other hand, we have
\[
\frac{1}{4}\varepsilon\sqrt{\varepsilon}> |\Im
z_{0} |= |{\sin}\psi | \biggl|R-\frac{1}{R}
\biggr|>\frac
{1}{2}\sqrt{\varepsilon}\frac{ |R^{2}-1 |}{R}.
\]
If $R>1$, then we deduce that $2R^{2}-\varepsilon R-2<0$. Therefore,
the number $R$ must be bounded from above by the positive
$x$-intercept of the parabola $Y=2X^{2}-\varepsilon X-2$. This means
that $R$ is no more than $1.026$ because $\varepsilon<0.1$.

Recall, from (\ref{eq:3.2}), that
\[
F_{n+1}(z)=z-G_{\nu_{n}\boxplus\gamma_{t}}(z),\qquad
z\in\mathbb{C}^{+}\cup\mathbb{R},
\]
where the sequence $\nu_{n}$
converges weakly to $\delta_{0}$ as $n\rightarrow\infty$. We will
derive an approximant for the function
$G_{\nu_{n}\boxplus\gamma_{t}}(z)$ as in (\ref{eq:3.7}) below. By
Proposition \ref{prop21} and Lemma \ref{lem33}, there exists $N=N(\varepsilon)>0$
such that, for $n\geq N$, the following conditions are satisfied:
\begin{enumerate}[(a)]
\item[(a)] $32/(n+1)<\varepsilon\sqrt{\varepsilon}$;
\item[(b)] the $R$-transform $R_{\nu_{n}}$ is defined in $D_{\theta,2}$
and the sequence $R_{\nu_{n}}$ converges uniformly on
$G_{\gamma} (K_{\varepsilon} )$ to zero as
$n\rightarrow\infty$;\vspace*{12pt}
\item[(c)]
\begin{eqnarray*}
\\[-43pt]
 |R_{\nu_{n}}(w) |<\frac{\varepsilon\sqrt{\varepsilon
}}{16},\qquad
w\in G_{\gamma} (K_{\varepsilon} );
\end{eqnarray*}
\item[(d)] the function $G_{\nu_{n}\boxplus\gamma_{t}}(z)$ has an
analytical
extension (which we still denote by $G_{\nu_{n}\boxplus\gamma_{t}}$)
to a neighborhood $\mathcal{V}$ of the interval
$[-2+2\varepsilon,2-2\varepsilon]$ and
$G_{\nu_{n}\boxplus\gamma_{t}}\neq0$ in $\mathcal{V}$.
\end{enumerate}
The additivity of $R$-transform shows that
\begin{eqnarray*}
G_{\nu_{n}\boxplus\gamma_{t}}^{-1}(w) & = & G_{\nu
_{n}}^{-1}(w)+G_{\gamma_{t}}^{-1}(w)-\frac{1}{w}\\
& = & w+\frac{1}{w}-\frac{1}{n+1}w+R_{\nu_{n}}(w),
\end{eqnarray*}
wherever the function $R_{\nu_{n}}$ is defined. Hence, (b) shows that
the function $G_{\nu_{n}\boxplus\gamma_{t}}^{-1}$ is defined in
$D_{\theta,2}$ for $n\geq N$. Moreover, after replacing $w$ by
$G_{\gamma}(z)$ and making use of (\ref{eq:2.1}), we conclude that
%
%
\begin{equation}\label{eq:3.5}
f_{n}(z)=G_{\nu_{n}\boxplus\gamma_{t}}^{-1} (G_{\gamma
}(z)
)=z+r_{n}(z),\qquad
z\in K_{\varepsilon},
\end{equation}
where the sequence
\[
r_{n}(z)=R_{\nu_{n}} (G_{\gamma}(z) )-\frac
{1}{n+1}G_{\gamma}(z)
\]
converges uniformly on $K_{\varepsilon}$ to zero as
$n\rightarrow\infty$ and the estimate
\[
 |r_{n}(z) |<\frac{\varepsilon\sqrt{\varepsilon}}{8}
\]
holds uniformly for $z\in K_{\varepsilon}$ and $n\geq N$.

The uniform bound of $r_{n}$ and (\ref{eq:3.5}) imply that the
rectangle $\widetilde{K}_{\varepsilon}=\{
x+iy\dvtx x\in[-2+2\varepsilon,2-2\varepsilon], |8y
|<\varepsilon
\sqrt{\varepsilon}\}$
is contained in the set $f_{n} (K_{\varepsilon} )$. Then
Rouch\'{e}'s theorem shows that each function $f_{n}$ has an
analytic inverse $f_{n}^{-1}$ defined in
$\widetilde{K}_{\varepsilon}$. Moreover, by inverting
(\ref{eq:3.5}), we have
\[
f_{n}^{-1}(z)=z+O ( |r_{n}(z) | )
\]
for every
$z\in\widetilde{K}_{\varepsilon}$ and $n\geq N$. Meanwhile, (d)
shows that the composition $G_{\gamma}^{-1}\circ
G_{\nu_{n}\boxplus\gamma_{t}}$ is defined and analytic in
$\mathcal{V}$, and hence it coincides with the function $f_{n}^{-1}$
on the interval $[-2+2\varepsilon,2-2\varepsilon]$. Therefore, we
conclude that there exist constants $Q>0$ and
$\delta=\delta(\varepsilon)\in(0,\varepsilon\sqrt{\varepsilon}/8)$,
both independent of $n$, such that
%
%
\begin{equation}\label{eq:3.6}
G_{\gamma}^{-1} (G_{\nu_{n}\boxplus\gamma_{t}}(z)
)=f_{n}^{-1}(z)=z+g_{n}(z)
\end{equation}
for every $z$ in the rectangle $K=\{
x+iy\dvtx x\in[-2+2\varepsilon,2-2\varepsilon], |y |<\delta\}$,
where $ |g_{n} |\leq Q |r_{n} |$ on $K$. Note that
the Cauchy estimate shows that
$ |G_{\gamma}^{\prime}(z) |\leq\varepsilon^{-1}$ on~$K$.
Therefore, by applying $G_{\gamma}$ on (\ref{eq:3.6}), we get
%
%
\begin{equation}\label{eq:3.7}
G_{\nu_{n}\boxplus\gamma_{t}}(z)=G_{\gamma} \bigl(z+g_{n}(z)
\bigr)=G_{\gamma}(z)+l_{n}(z),\qquad
z\in K,  n\geq N,
\end{equation}
where
$ |l_{n}(z) |\leq {Q\sup_{z\in
K}} |G_{\gamma}^{\prime}(z) | |r_{n}(z) |\leq
Q\varepsilon^{-1} |r_{n}(z) |$ on $K$.

By (\ref{eq:3.2}), (\ref{eq:3.7}) and (\ref{eq:2.1}), the function
$G_{n+1}=1/F_{n+1}$ can be continued meromorphically to the set $K$. Moreover,
if $h_{n+1}$ denotes this continuation, then we have
\[
h_{n+1}(z)=\frac{G_{\gamma}(z)}{1-l_{n}(z)G_{\gamma}(z)},\qquad  z\in
K,  n\geq N.
\]
Therefore, after letting $n$ be large enough so that
$ |l_{n}(z)G_{\gamma}(z) |<1$, we conclude that the
function $h_{n+1}$ is the analytic continuation of $G_{n+1}$ to $K$,
and the uniform convergence property in the statement (ii) follows
because $l_{n}\rightarrow0$ uniformly on $K$ as
$n\rightarrow\infty$. The proof of (ii) is complete.

We next prove (iii). The proof of (ii) shows that, for large $n$, we
have
\[
 |G_{n+1}(x)-G_{\gamma}(x) |<\varepsilon,\qquad
x\in[-2+2\varepsilon,2-2\varepsilon].
\]
Hence, by the inversion
formula, the density function $d\mu_{n+1}/dx$ is uniformly close to
$d\gamma/dx$ on $[-2+2\varepsilon,2-2\varepsilon]$ for sufficiently
large $n$. Combining with (\ref{eq:3.2}) and~(\ref{eq:2.3}), we then
have
\[
\Im G_{n+1}(x)=\frac{\Im
G_{\nu_{n}\boxplus\gamma_{t}}(x)}{ [x-\Re
G_{\nu_{n}\boxplus\gamma_{t}}(x) ]^{2}+ [\Im
G_{\nu_{n}\boxplus\gamma_{t}}(x) ]^{2}},\qquad
x\in\mathbb{R},
\]
for large $n$. Hence, we see that the function
$d\mu_{n+1}/dx$ is continuous for such $n$.

In order to finish the proof, we need to estimate $d\mu_{n+1}/dx$ on
the set $\mathbb{R}\setminus[-2+2\varepsilon,2-2\varepsilon]$. Let
us take $x=\psi_{t}(u)$ and $v_{t}(u)=-t\Im
G_{\nu_{n}\boxplus\gamma_{t}}(x)$ as in Theorem \ref{theo22}. Then Lemma \ref{lem32}
and (\ref{eq:2.3}) imply that $ |u |>1-2\varepsilon$ and
\[
 |x-\Re
G_{\nu_{n}\boxplus\gamma_{t}}(x) |\geq |x |-\tfrac
{7}{6}, \qquad |x |>2-2\varepsilon,
n\geq3.
\]
Hence, we obtain the bound:
%
%
\begin{equation}\label{eq:3.8}
|\Im G_{n+1}(x) |\leq\frac{72v_{t}(u)}{( |x
|+1)^{2}},\qquad
 |u |>1-2\varepsilon,  |x |>2-2\varepsilon,
\end{equation}
for every $n\geq3$.

By the inversion formula, it remains to show that $v_{t}$ is
uniformly small on $\{ u\dvtx |u |>1-2\varepsilon\}$ as
$n\rightarrow\infty$. This can be done as follows. For each
$n\geq1$, we introduce the probability tail-sum
$k_{n}=1-\nu_{n}((-\varepsilon,\varepsilon))$. From (\ref{eq:2.2}),
we have, for $v_{t}(u)>0$ and $ |u |>1-2\varepsilon$,
that
\begin{eqnarray*}
\frac{1}{t} & = & \int_{-\infty}^{\infty}\frac
{1}{(u-x)^{2}+v_{t}(u)^{2}} \, d\nu_{n}(x)\\
& \leq& \int_{ |x |<\varepsilon}\frac{1}{(1-3\varepsilon
)^{2}+v_{t}(u)^{2}} \, d\nu_{n}(x)+\int_{ |x |\geq
\varepsilon
}\frac{1}{v_{t}(u)^{2}} \, d\nu_{n}(x)\\
& = & \frac{1-k_{n}}{(1-3\varepsilon)^{2}+v_{t}(u)^{2}}+\frac
{k_{n}}{v_{t}(u)^{2}}.
\end{eqnarray*}
Therefore, we see that the number $v_{t}(u)^{2}$ is bounded above by
the positive $x$-intercept of the parabola $Y=X^{2}-bX-c$, where
$b=t-(1-3\varepsilon)^{2}$ and $c=tk_{n}(1-3\varepsilon)^{2}\geq0$.
In other words, we have
\[
v_{t}(u)^{2}\leq\frac{b+\sqrt{b^{2}+4c}}{2}\leq |b
|+\sqrt
{k_{n}},\qquad  |u |\geq1-2\varepsilon,
n\geq1.
\]
Hence, the desired conclusion for $v_{t}$ follows from the
facts that $\lim_{n\rightarrow\infty}k_{n}=0$ and
$\lim_{n\rightarrow\infty}b=6\varepsilon-9\varepsilon^{2}$. Since
$\varepsilon$ is arbitrary, this completes the proof of (iii).
\end{pf}
\begin{rem*}
Observe first that, in Theorem \ref{theo34}, the function
$(\overline{h_{n}(\overline{z})}-h_{n}(z))/ 2\pi i$ is a complex
analytic extension of $d\mu_{n}/dx$ to $K$ for $n\geq N$. Hence, for
any $\varepsilon>0$ and $k\geq1$, Theorem \ref{theo34} actually implies that
the $k$th derivatives $d^{k}\mu_{n}/dx^{k}$ converge to
$d^{k}\gamma/dx^{k}$ uniformly on $[-2+\varepsilon,2-\varepsilon]$
as $n\rightarrow\infty$.

Secondly, if the measure $\mu$ has a bounded support, then
so does the measure~$\mu_{n}$. In this case, the equation
(\ref{eq:3.2}) shows that the support of the measure $\mu_{n}$ can
not deviate too much from the interval $[-2,2]$. This property has
been noted in \cite{Superconvergence}, where it was shown that, for
any $\varepsilon>0$ there exists $N>0$ such that the support of
$\mu_{n}$ is contained in $[-2-\varepsilon,2+\varepsilon]$ for
$n\geq N$. An earlier result of Voiculescu \cite{Voiadd},
Lemmas \ref{lem31} and \ref{lem32}, provides a precise estimate for the support of a free
convolution of compactly supported measures (see also
\cite{Kargin}). In particular, Voiculescu's result shows that the
support of $\mu_{n}$ is contained in the interval
$[-2-(L/\sqrt{n}),2+(L/\sqrt{n})]$ for $n\geq1$, where
$L=\sup\{ |x |\dvtx x\in\operatorname{supp}(\mu)\}$.
\end{rem*}

We conclude this section with a global form of the central limit
theorem.
\begin{theorem}\label{theo35}
Let $\mu$ and $\mu_{n}$ be as in Theorem \ref{theo34}. Then we have
\[
\lim_{n\rightarrow\infty} \biggl\Vert
\frac{d\mu_{n}}{dx}-\frac{d\gamma}{dx} \biggr\Vert
_{L^{p} (\mathbb{R} )}=0
\]
for every $p>1/2$.
\end{theorem}
\begin{pf}
Fix $p>1/2$. Let $\varepsilon\in(0,1/10)$ be given. From the
inversion formula, we have
\begin{eqnarray*}
\biggl\Vert\frac{d\mu_{n+1}}{dx}-\frac{d\gamma}{dx} \biggr\Vert
_{L^{p} (\mathbb{R} )}^{p} & = & \frac{1}{\pi^{p}}\int
_{-2+2\varepsilon}^{2-2\varepsilon} |\Im G_{n+1}(x)-\Im
G_{\gamma
}(x) |^{p} \, dx\\
& & {} +\frac{1}{\pi^{p}}\int_{ |x |>2-2\varepsilon}
|\Im
G_{n+1}(x)-\Im G_{\gamma}(x) |^{p} \, dx.
\end{eqnarray*}
Then Theorem \ref{theo34}(ii) shows that within the interval
$[-2+2\varepsilon,2-2\varepsilon]$ the integrand tends uniformly to
zero as $n\rightarrow\infty$ and so the contribution of
$[-2+2\varepsilon,2-2\varepsilon]$ tends to zero. As to the second
integral,
\[
I_{n}=\int_{ |x |>2-2\varepsilon} |\Im
G_{n+1}(x)-\Im G_{\gamma}(x) |^{p} \, dx,
\]
(\ref{eq:3.8})
implies, for $n\geq3$, that
\begin{eqnarray*}
I_{n} & \leq& \int_{ |x |>2-2\varepsilon}2^{p}
\bigl( |\Im
G_{n+1}(x) |^{p}+ |\Im G_{\gamma}(x) |^{p} \bigr)\,
dx\\
& \leq& 144^{p} \bigl(\sqrt{ |b |}+\sqrt[4]{k_{n}}
\bigr)^{p}\int_{ |x |>2-2\varepsilon}\frac{1}{( |x
|+1)^{2p}} \, dx\\
& &{} +\int_{2\geq |x |>2-2\varepsilon} (4-x^{2}
)^{{p}/{2}} \, dx,
\end{eqnarray*}
where $b$ and $k_{n}$ are as in the proof of Theorem \ref{theo34}. Then the
proof is completed by taking $n\rightarrow\infty$, then
$\varepsilon\rightarrow0^{+}$.
\end{pf}

\section{Entropic central limit theorem}

In this section, we study the convergence in the free central limit
theorem in terms of entropy. Our motivation is a result of Barron,
which we now review as follows.

The classical entropy of a measure $\rho\in\mathcal{M}$ with density
$f$ is defined as $H(\rho)=-\int_{-\infty}^{\infty}f(x)\log f(x)
dx$, provided the positive part of the integral is finite. Thus we
have $H(\rho)\in[-\infty,\infty)$. It is well known that the
standard Gaussian distribution $G$ has the largest entropy among all
probability measures on $\mathbb{R}$ with variance one. Let $\mu$,
$D_{1/\sqrt{n}}\mu$ and $\mu_{n}$ be as in Section \ref{sec3}, and further
let $\rho_{n}$ be the classical convolution
\[
\rho_{n}=\underbrace{D_{{1}/{\sqrt{n}}}\mu*D_{{1}/{\sqrt
{n}}}\mu
*\cdots*D_{{1}/{\sqrt{n}}}\mu}_{n\ \mathrm{times}}.
\]
Barron \cite{Barron} showed that if the entropy $H(\rho_{n})$ ever
becomes different from $-\infty$, then the sequence $H(\rho_{n})$
converges to $H(G)=(1/2)\log2\pi e$. A detailed exposition of the
connections between the entropy and central limit theorems can be
found in \cite{Johnson}.

Let $\nu$ be a probability measure on $\mathbb{R}$. The quantity
\[
\chi(\nu)=\int_{-\infty}^{\infty}\int_{-\infty}^{\infty}{\log}
|x-y |\,
d\nu(x)\,d\nu(y)+\frac{3}{4}+\frac{1}{2}\log2\pi,
\]
called \textit{free
entropy}, was discovered by Voiculescu in \cite{VoiI} to be a good
substitute for entropy in free probability theory. Remarkably, free
entropy $\chi$ behaves like the classical entropy $H$ in many
instances. For example, the free entropy is maximized by the
standard semicircular law $\gamma$ [with the value
$\chi(\gamma)=(1/2)\log2\pi e$] among all probability measures with
variance one \cite{VoiIV,HPMaxEntropy} and both $\chi(\mu_{n})$ and
$H(\rho_{n})$ are monotonically increasing
\cite{Shlya,Schultz,ABBN}. We also refer to \cite{VoiSurvey} for a
survey of free entropy. Our goal here is to prove a free analogue of
Barron's result.

The key ingredient in the proof of our result is a free logarithmic
Sobolev inequality first proved by Voiculescu \cite{VoiV}, which
requires a further assumption on the regularity of the measure
$\nu$. Suppose that the measure $\nu$ has a density $p$ in
$L^{3} (\mathbb{R} )$, and let us assume the support of
$\nu$ is bounded for the moment. Then the \textit{free Fisher
information} of the measure $\nu$ is
\[
\Phi(\nu)=\frac{4\pi^{2}}{3}\int_{-\infty}^{\infty}p(x)^{3} \, dx;
\]
see \cite{VoiV}, Proposition 7.9, where the following inequality is
proved:
%
%
\begin{equation}\label{eq:4.1}
\chi(\nu)\geq\frac{1}{2}\log \biggl(\frac{2\pi
e}{\Phi(\nu)} \biggr).
\end{equation}
Using random
matrix approximations, the inequality (\ref{eq:4.1}) is further
shown to hold for measures with finite variance, not necessarily
having a bounded support (see~\cite{BianeLogSobolev}). Thus, we can
use this generalized version of (\ref{eq:4.1}) to obtain the
following:
\begin{theorem}\label{theo41}
Let $\mu$ and $\mu_{n}$ be as in Theorem \ref{theo34}. Then the free
entropy $\chi(\mu_{n})$ converges to the semicircular entropy
$(1/2)\log2\pi e$.
\end{theorem}
\begin{pf}
Notice first that $\Phi(\gamma)=1$. Theorem \ref{theo35} shows that the free
Fisher information $\Phi(\mu_{n})$ converges to $\Phi(\gamma)=1$ as
$n\rightarrow\infty$. Since $\mu_{n}$ has variance one, we conclude,
from (\ref{eq:4.1}), that $\chi(\mu_{n})>-\infty$ and
\[
\chi(\gamma)\geq\chi(\mu_{n})\geq\chi(\gamma)-\tfrac{1}{2}\log
 (\Phi
(\mu_{n}) )
\]
for sufficiently large $n$. The result follows immediately.
\end{pf}

Finally, we conclude this paper with the following remark. A version
of Theorem \ref{theo41} for bounded variables satisfying more restrictive
conditions was noted in~\cite{Johnson}. Here, we only require the
variance to be one, which is the most general condition to have the
free central limit theorem (see \cite{PataCLT}). Actually, Theorem
\ref{theo41} can be extended to any probability measure with finite nonzero
variance via a standard translating and scaling procedure by
observing that the free entropy $\chi$ is invariant under the
translation and
\[
\chi (D_{a}\nu )=\chi (\nu )+\log a, \qquad a>0.
\]

\section*{Acknowledgment}
The author would like to thank Professors
Hari\break Bercovici and Roland Speicher for their encouragement and
valuable help during the course of the investigation. He is also
grateful to Professor Michael Anshelevich for a reference, and to a
referee for helpful comments on the preliminary version of this
work.

%

%
\printaddresses

\end{document}